# CRAMÉR'S ESTIMATE FOR A REFLECTED LÉVY PROCESS[1]

By R. A. Doney and R. A. Maller

*University of Manchester and Australian National University*

The natural analogue for a Lévy process of Cramér's estimate for a reflected random walk is a statement about the exponential rate of decay of the tail of the characteristic measure of the height of an excursion above the minimum. We establish this estimate for any Lévy process with finite negative mean which satisfies Cramér's condition, and give an explicit formula for the limiting constant. Just as in the random walk case, this leads to a Poisson limit theorem for the number of "high excursions."

**1. Introduction.** The reflected process $R = (R_n,\ n \geq 0)$ formed from a random walk $S = (S_n, n \geq 0)$ by setting

$$R_n = S_n - I_n \qquad \text{where } I_n = \min_{i \leq n} S_i,\ n \geq 0,$$

arises in many areas of applied probability, including queuing theory, risk theory and mathematical genetics; see, for example, [1, 2, 6] and [8]. In these applications, the current maximum of $R$,

$$M_n^{(R)} = \max_{j \leq n} R_j = \max_{0 \leq i \leq j \leq n} \{S_j - S_i\},$$

is of particular interest; in the genetics context this is called the maximal segmental score. The segments here correspond to excursions of the random walk above its minimum, and, more generally, the behavior of

$$N(y, n) = \#\{\text{excursions completed by time } n \text{ whose heights exceed } y\}$$

is important. A classical result in [7] asserts that if the random walk has a finite negative drift and Cramér's condition is satisfied, then $N(y, n)$ has a

Received May 2003; revised October 2004.
[1]Supported in part by ARC Grant DP0210572.
*AMS 2000 subject classifications.* 60G51, 60G17.
*Key words and phrases.* Maximum of reflected process, maximal segmental score, Poisson limit theorem, high excursions.

limiting Poisson distribution when $n, y \to \infty$ in such a way that $ne^{-y}$ converges to a positive constant. In [8] extensions of this result to processes other than random walks were given, including compound Poisson processes and Brownian motion with negative drift. We show here that these extensions actually hold when the underlying process is any Lévy process with finite negative drift which satisfies Cramér's condition.

**2. The random walk case.** Let $S$ be any random walk with finite negative mean which satisfies Cramér's condition, namely $E(e^{\gamma S_1}) = 1$ for some $\gamma \in (0, \infty)$. Let $h_i, i = 1, \ldots,$ denote the height of the $i$th excursion above the minimum, that is,

$$h_i = \max_{0 \le n \le T_i - T_{i-1}} \{S_{T_{i-1}+n} - S_{T_{i-1}}\},$$

where $T_i$ is the $i$th strict descending ladder time, with $T_0 = 0$, and let $M_\infty$ denote the all-time maximum of the walk. Also write $T_i^+$, $H_i^+ = S(T_i^+)$ for the weak increasing ladder times and heights and $H_i = |S(T_i)|$ for the decreasing ladder heights. Then Cramér's famous estimate states that

$$\lim_{x\to\infty} e^{\gamma x} P(M_\infty > x) = C = \frac{P(H_1^+ < \infty)}{\gamma E(H_1^+ e^{\gamma H_1^+}; H_1^+ < \infty)}. \tag{1}$$

(See, e.g., page 413 of [6].) Obviously, the $h_i$ are independent, identically distributed and it is easy to deduce from the above that, in the nonlattice case,

$$e^{\gamma x} P(h_1 > x) \to K := C\{1 - E(e^{-\gamma H_1})\}. \tag{2}$$

To see this, observe the identity

$$P(M_\infty > x) = P(h_1 > x) + \int_0^\infty P(h_1 \le x, H_1 \in dy) P(\widetilde{M}_\infty > x + y), \tag{3}$$

where $\widetilde{M}_\infty$ is an independent copy of $M_\infty$, multiply by $e^{\gamma x}$, and let $x \to \infty$; this argument is due to Iglehart [7].

Now introduce the finite constant $\alpha = ET_1$, so that the strong law implies that

$$\frac{\sup\{i : T_i \le n\}}{n} \overset{\text{a.s.}}{\to} \frac{1}{\alpha} \qquad \text{as } n \to \infty.$$

Using (2), it is then easy to see that if $n$ and and $y \to \infty$ in such a way that $ne^{-\gamma y} \to \alpha\lambda/K$, then the number of the $h_i$ which exceed $y$ and occur by time $n$ has a limiting Poisson($\lambda$) distribution.

**3. The Lévy process case.** Now let $X$ be any Lévy process with $EX_1 \in (-\infty, 0)$ which satisfies Cramér's condition, namely

$$E(e^{\gamma X_1}) = 1 \qquad \text{for some } \gamma \in (0, \infty),$$

and define $Y$, the process reflected in its infimum, by

$$Y_s = X_s - I_s \qquad \text{where } I_s = \inf_{0 \le u \le s} X_u,\ s \ge 0.$$

Let $L = (L_t, t \ge 0)$ and $L^{-1} = (L_t^{-1}, t \ge 0)$ denote the local time process of $Y$ at 0 and its right continuous inverse, and write

$$\varepsilon_t(u) = Y(L_t^{-1} + u) - Y(L_t^{-1}), \qquad u \ge 0, \quad \text{and} \quad \xi_t = \inf\{u : \varepsilon_t(u) \le 0\}.$$

Then $\varepsilon_t = (\varepsilon_t(u), 0 \le u < \xi_t)$ is the excursion above the minimum at local time $t$, and

$$h_t = \sup(\varepsilon_t(u), 0 \le u < \xi_t)$$

is the height of this excursion. Thus, the Poisson point process $(h_t, t \ge 0)$ is the continuous time version of $(h_i, i \ge 1)$ for a random walk. If $\eta$ denotes the characteristic measure of $(h_t, t \ge 0)$, then the statement which corresponds to (2) in this context is that

$$e^{\gamma x} \eta((x, \infty)) \to K^* \qquad \text{as } x \to \infty. \tag{4}$$

Our aim is to establish (4), determine the constant $K^*$ and deduce that the number of excursions of $Y$ away from 0 with heights exceeding $y$ which take place by time $t$ satisfies a Poisson limit theorem. Of course, the reason that this is potentially more difficult than the random walk case is that there may be an infinite number of excursions in any finite time interval.

The starting point, naturally, is the following known Lévy process version of (1); (see [4]): with $S_\infty = \sup_{0 \le t < \infty} X_t$ denoting the all-time supremum and $H^+ = (H_t^+, t \ge 0)$ the increasing ladder-height process,

$$e^{\gamma x} P(S_\infty > x) \to C^* = \frac{\beta}{\gamma m} \qquad \text{as } x \to \infty, \tag{5}$$

where

$$\beta = -\log P(H_1^+ < \infty) \quad \text{and} \quad m = E(H_1^+ e^{\gamma H_1^+}; H_1^+ < \infty). \tag{6}$$

Actually, it is easy, by applying the random walk argument to the independent, identically distributed sequence of excursions heights which exceed some fixed $\delta > 0$, to deduce from (5) that (4) holds, but with a value of $K^*$ which apparently depends on $\delta$. Thus, the proof of the following theorem, which is quite delicate, is essentially a matter of justifying an interchange of limits. The inverse local time process $L^{-1} = (L_t^{-1}, t \ge 0)$ is the Lévy version of the descending ladder time process, and the corresponding ladder height process is defined by

$$H_t = |X(L_t^{-1})| = |I(L_t^{-1})| \qquad \text{where } I_t = \inf_{0\le s\le t} X_s.$$

THEOREM 1. (i) *Let* $\phi(\theta) = -\log E(e^{-\theta H_1})$ *denote the exponent of the subordinator* $H$. *Then, as* $x\to\infty$,

$$e^{\gamma x}\eta((x,\infty)) \to K^* := \phi(\gamma)C^* = \frac{\phi(\gamma)\beta}{\gamma m}.$$

(ii) *Introduce the finite constant* $\alpha^* = EL_1^{-1}$, *and let* $N(y,t)$ *denote the number of excursions of* $Y$ *with heights exceeding* $y$ *which occur by time* $t$. *Let* $t,y\to\infty$ *in such a way that* $te^{-\gamma y}\to \alpha^*\lambda/K^*$. *Then* $N(y,t)$ *has a limiting* Poisson$(\lambda)$ *distribution.*

PROOF. Write $S_t = \sup_{0\le u\le t} X_u$ and $\widehat{S}_t = S(L_t^{-1})$. Then applying the Markov property at the stopping time $L_t^{-1}$, we see that

$$P(S_\infty \le x) = P\{(\widehat{S}_t \le x)\cap(\widetilde{S}_\infty \le x + H_t)\},$$

where $\widetilde{S}_\infty$ is independent of $\widehat{S}_t$ and $H_t$ and has the distribution of $S_\infty$. Note, however, that $H_t$ and $\widehat{S}_t$ are dependent. Subtracting the term $P(\widehat{S}_t \le x)P(S_\infty\le x)$ gives

$$P(\widehat{S}_t > x)P(S_\infty \le x) = P\{(\widehat{S}_t \le x)\cap(x < \widetilde{S}_\infty \le x + H_t)\}. \tag{7}$$

Our first step is to examine the behavior of $e^{\gamma x}$ $P(\widehat{S}_t > x)$ as $x\to\infty$ for fixed $t$. Note first that in view of (5) and the fact that $P(\widehat{S}_t > x)\to 0$ as $x\to\infty$ for fixed $t$,

$$e^{\gamma x}P\{(\widehat{S}_t > x)\cap(x<\widetilde{S}_\infty\le x+H_t)\} \le e^{\gamma x}P(\widehat{S}_t > x)P(\widetilde{S}_\infty > x)\to 0.$$

Thus, as $x\to\infty$, for each fixed $t$,

$$\begin{aligned}
&P(x<\widetilde{S}_\infty\le x+H_t)\\
&\qquad = e^{\gamma x}P\{(\widehat{S}_t > x)\cap(x<\widetilde{S}_\infty\le x+H_t)\} + o(1) && (8)\\
&\qquad = e^{\gamma x}P(\widehat{S}_t > x)P(S_\infty\le x) + o(1) && (9)\\
&\qquad = e^{\gamma x}P(\widehat{S}_t > x) + e^{\gamma x}P(\widehat{S}_t > x)P(S_\infty > x) + o(1) && (10)\\
&\qquad = e^{\gamma x}P(\widehat{S}_t > x) + o(1). && (11)
\end{aligned}$$

Next recall from [4] that we can write $e^{\gamma x}P(S_\infty\in dx) = \beta U(dx)$, where $U$ is the renewal measure corresponding to a distribution on $[0,\infty)$ with the finite mean $m$ given in (6). So we can write

$$e^{\gamma x}P(x<\widetilde{S}_\infty\le x+H_t) = \beta E\left(\int_0^{H_t} e^{-\gamma y}U(x+dy)\right) := \beta E Z_t(x). \tag{12}$$

Since $mU(x+dy)$ converges weakly to a Lebesgue measure, we see that, for each fixed $\omega$ and $t$,

$$Z_t(x,\omega) \to m^{-1}\int_0^{H_t(\omega)} e^{-\gamma y}\,dy = \frac{1-e^{-\gamma H_t(\omega)}}{m\gamma} \qquad \text{as } x\to\infty. \tag{13}$$

Also, using the subadditivity of $U$ and Erickson's bounds (see [5]), we get

$$Z_t(x,\omega) \le U(x+H_t(\omega)) - U(x) \le U(H_t(\omega)) \le \frac{2H_t(\omega)}{m}.$$

Since $EX_1 \in (-\infty, 0)$, we have $EH_t < \infty$ for any fixed $t$, so by dominated convergence, we conclude from (8), (12) and (13) that

$$\begin{aligned}\lim_{x\to\infty} e^{\gamma x} P(\widehat{S}_t > x) &= \beta E \lim_{x\to\infty} Z_t(x) \\ &= \frac{\beta E(1-e^{-\gamma H_t})}{m\gamma} = \frac{\beta(1-e^{-t\phi(\gamma)})}{m\gamma}.\end{aligned} \tag{14}$$

To connect this to the asymptotic behavior of $e^{\gamma x}\overline{\eta}(x)$, where $\overline{\eta}(x) = \eta((x,\infty))$, we need the following observation; the event $\widehat{S}_t \le x$ occurs if and only if each excursion at local time $s$, for all $s \le t$, has height $\le H_s + x$. Hence, by a standard application of the compensation formula (see [3], page 7), we have

$$t^{-1}P(\widehat{S}_t > x) = t^{-1}E\biggl(1-\exp\biggl\{-\int_0^t \overline{\eta}(x+H_s)\,ds\biggr\}\biggr). \tag{15}$$

Using the bound $1-e^{-x} \le x$ in this shows that, for any $t>0$,

$$t^{-1}e^{\gamma x}P(\widehat{S}_t > x) \le t^{-1}e^{\gamma x}E\int_0^t \overline{\eta}(x+H_s)\,ds \le e^{\gamma x}\overline{\eta}(x).$$

From (14) we see that

$$\liminf_{x\to\infty} e^{\gamma x}\overline{\eta}(x) \ge \lim_{t\downarrow 0}\lim_{x\to\infty} t^{-1}e^{\gamma x}P(\widehat{S}_t > x) = \lim_{t\downarrow 0}\frac{\beta(1-e^{-t\phi(\gamma)})}{m\gamma t} = K^*.$$

Also, using the bound $1-e^{-x} \ge x - x^2/2$ in (15) gives, for any fixed $\varepsilon > 0$ and all $t>0$,

$$\begin{aligned} &t^{-1}e^{\gamma x}P(\widehat{S}_t > x) \\ &\quad \ge t^{-1}e^{\gamma x}\biggl(E\int_0^t \overline{\eta}(x+H_s)\,ds - \frac{1}{2}E\biggl\{\int_0^t \overline{\eta}(x+H_s)\,ds\biggr\}^2\biggr) \\ &\quad \ge e^{\gamma x}\biggl(E\overline{\eta}(x+H_t) - \frac{t}{2}\overline{\eta}(x)^2\biggr) \\ &\quad \ge e^{\gamma x}\overline{\eta}(x+\varepsilon)P(H_t \le \varepsilon) - \frac{t}{2}e^{\gamma x}\overline{\eta}(x)^2.\end{aligned}$$

Rearranging, using (14), the fact that $P(H_t \le \varepsilon) \to 1$ as $t \downarrow 0$, and noting that the final term above is $o(\limsup_{x\to\infty} e^{\gamma x}\overline{\eta}(x))$, we see that

$$\begin{aligned}\limsup_{x\to\infty} e^{\gamma x}\overline{\eta}(x) &= e^{\gamma\varepsilon} \limsup_{x\to\infty} e^{\gamma x}\overline{\eta}(x+\varepsilon)\\ &\le e^{\gamma\varepsilon} \lim_{t\downarrow 0} \lim_{x\to\infty} e^{\gamma x} t^{-1} P(\widehat{S}_t > x)/P(H_t \le \varepsilon)\\ &= e^{\gamma\varepsilon} K^*.\end{aligned}$$

Since $\varepsilon$ is arbitrary, the proof of (i) is complete.

For (ii), note that standard properties of Poisson point processes show that $N(y,t)$ has a Poisson distribution with parameter $(L_t\overline{\eta}(y))$. It follows from the strong law for subordinators ([3], page 92) that $L_t \sim t/\alpha^*$ as $t \to \infty$, so the conclusion follows. □

Remark 2. It is easy to adapt our arguments to extend the convergence in (ii) of Theorem 1 to joint convergence in law for different values of $\lambda$, and this leads to a process version of our result. (We are grateful to a referee for this observation.)

## REFERENCES

[1] Asmussen, S. (1982). Conditioned limit theorems relating a random walk to its associate, with applications to risk reserve processes and the $GI/GI/1$ queue. *Adv. in Appl. Probab.* **14** 143–170. MR644012

[2] Asmussen, S. (1982). *Applied Probability and Queueing*. Wiley, New York.

[3] Bertoin, J. (1996). *Lévy Processes*. Cambridge Univ. Press. MR1406564

[4] Bertoin, J. and Doney, R. A. (1994). Cramér's estimate for Lévy processes. *Statist. Probab. Lett.* **21** 363–365. MR1325211

[5] Erickson, K. B. (1973). The strong law of large numbers when the mean is undefined. *Trans. Amer. Math Soc.* **185** 371–381. MR336806

[6] Feller, W. E. (1971). *An Introduction to Probability Theory and Its Applications* **2**, 2nd ed. Wiley, New York.

[7] Iglehart, D. L. (1972). Extreme values in the $GI/G/1$ queue. *Ann. Math. Statist.* **43** 627–635. MR305498

[8] Karlin, S. and Dembo, A. (1992). Limit distributions of maximal segmental score among Markov-dependent partial sums. *Adv. in Appl. Probab.* **24** 113–140. MR1146522

Department of Mathematics
University of Manchester
Oxford Road
Manchester M13 9PL
United Kingdom
e-mail: rad@ma.man.ac.uk

Centre for Mathematical Analysis
and School of Finance & Applied Statistics
Australian National University
Canberra, ACT
Australia
e-mail: Ross.Maller@anu.edu.au